\documentclass[10pt,notitlepage,twoside,a4paper]{amsart}
 \usepackage{amsfonts}

\usepackage{amsmath,amssymb,enumerate}

\usepackage{epsfig,fancyhdr,color}

\usepackage{amssymb}
\usepackage{amsmath,amsthm}
\usepackage{latexsym}
\usepackage{amscd}
\usepackage{psfrag}
\usepackage{graphicx}
\usepackage[latin1]{inputenc}
\usepackage[all]{xy}
\usepackage[mathcal]{eucal}

\definecolor{NoteColor}{rgb}{1,0,0}

\renewcommand{\textsc}{\textcolor{red}}

\newtheorem{theorem}{\rm\bf Theorem}[section]
\newtheorem{proposition}[theorem]{\rm\bf Proposition}
\newtheorem{lemma}[theorem]{\rm\bf Lemma}
\newtheorem{corollary}[theorem]{\rm\bf Corollary}
\newtheorem*{theorem 1}{\rm\bf Proposition 1}
\newtheorem*{theorem 2}{\rm\bf Proposition 2}

\theoremstyle{definition}

\theoremstyle{remark}

\newtheorem{example}[theorem]{\rm\bf Example}

\def\interieur#1{\mathord{\mathop{\kern 0pt #1}\limits^\circ}}

\title[Local comparison]{On local comparison between various metrics on Teichm\"uller spaces}

\author{D. Alessandrini}
\address{Daniele Alessandrini, Max-Planck-Institut f\"ur Mathematik, Vivatsgasse 7, 53111 Bonn, Germany}
\email{daniele.alessandrini@gmail.com}

\author{L. Liu}
\address{Lixin Liu, Department of Mathematics, Sun Yat-sen University, 510275, Guangzhou, P. R. China}
\email{mcsllx@mail.sysu.edu.cn}

\author{A. Papadopoulos}
\address{Athanase Papadopoulos, Institut de Recherche Math{\'e}matique Avanc\'ee,
Universit{\'e} de Strasbourg and CNRS,
7 rue Ren\'e Descartes,
 67084 Strasbourg Cedex, France} \email{athanase.papadopoulos@math.unistra.fr}
\date{\today}

\author{W. Su}
\address{Weixu Su, Department of Mathematics, Sun Yat-sen University, 510275, Guangzhou, P. R. China}
\email{suweixu@gmail.com}

\begin{document}

\begin{abstract} There are several Teichm\"uller spaces associated to a surface of infinite topological type, after the choice of a particular basepoint ( a complex or a hyperbolic structure on the surface). These spaces include the quasiconformal Teichm\"uller space, the length spectrum Teichm\"uller space, the Fenchel-Nielsen Teichm\"uller space, and there are others. In general, these spaces are set-theoretically different. An important question is therefore  to understand relations between these spaces. Each of these spaces is equipped with its own metric, and under some hypotheses, there are inclusions between these spaces. In this paper, we obtain local metric comparison results on these inclusions, namely, we show that the inclusions are locally bi-Lipschitz under certain hypotheses. To obtain these results, we use some hyperbolic geometry  estimates that give new results also for surfaces of finite type. We recall that in the case of a surface of finite type, all these Teichm\"uller spaces coincide setwise. In the case of a surface of finite type with no boundary components (and possibly with punctures), we show that the restriction of the identity map to any thick part of Teichm\"uller space is globally bi-Lipschitz with respect to the length spectrum metric and the classical Teichm\"uller metric on the domain and on the range respectively. In the case of a surface of finite type with punctures and boundary components, there is a metric on the Teichm\"uller space which we call the arc metric, whose definition is analogous to the length spectrum metric, but which uses lengths of geodesic arcs instead of lengths of closed geodesics. We show that the restriction of the identity map restricted to any ``relative thick" part of Teichm\"uller space is globally bi-Lipschitz, with respect to any of the three metrics: the length spectrum metric, the Teichm\"uller metric and the arc metric on the domain and on the range.

 \bigskip

\medskip

 \bigskip

 \noindent Liu and Su are partially supported by NSFC grant No. 10871211.

 \end{abstract}

\maketitle

\tableofcontents

\section{Introduction}\label{intro}

The paper concerns surfaces of finite and of infinite topological type. The results are not identical for both cases, and we treat them separately. We start with the case of surfaces of infinite type.

Let $\Sigma$ be a surface of infinite topological type, that is, a surface obtained by gluing a countably  infinite number of generalized pairs of pants along their boundary components. Here, a generalized pair of pants is a sphere with three holes, a hole being either a point removed or an open disk removed.
There are several Teichm\"uller spaces associated with such a surface $\Sigma$, with several inclusions between these spaces, and different metrics on them. We are interested in comparing these metrics, in cases where a comparison can be done. This paper is a continuation of the work done in  \cite{ALPSS} and \cite{completeness}, in which we introduced a space we called the Fenchel-Nielsen Teichm\"uller space, which is equipped with a metric we called the Fenchel-Nielsen metric. We compared this metric with the Teichm\"uller metric.  The definition of the Fenchel-Nielsen Teichm\"uller space of $\Sigma$ depends on the choice of a basepoint for this space and of a pair of pants decomposition of the surface.
 Our work is also in the spirit of \cite{LP}, in which we studied the various metrics on Teichm\"uller spaces of surfaces of infinite topological type. In the present paper, we mainly consider the question of local metric comparison between the Fenchel-Nielsen metric, the quasiconformal metric and the length spectrum metric.

Teichm\"uller spaces can be seen as parameter spaces for conformal structures on $\Sigma$. We will always consider these conformal structures as endowed with their \emph{intrinsic metric}. This is a hyperbolic metric in the given conformal class, and it  was defined by Bers. In the case of a surface with no boundary components (that is, no ends whose neighborhoods are conformally equivalent to annuli) but which may have punctures (that is, ends with neighborhoods conformally to punctured discs), the intrinsic metric coincides with the Poincar\'e metric. But in the case of a surface with boundary components the two metrics do not coincide,
 see the end of Section 4 of \cite{ALPSS} for the definition and a discussion. Endowing the Riemann surface with its intrinsic metric will allow us to use techniques from hyperbolic geometry, like the existence of hyperbolic pair of pants decompositions and of Fenchel-Nielsen coordinates. We note that in order to define Fenchel-Nielsen coordinates, we need to show that given a topological pair of pants decomposition $\mathcal{P}=\{C_i\}_{i=1,2,\ldots}$ of $\Sigma$ and a hyperbolic metric on $\Sigma$, there exists a unique geodesic pair of pants decomposition in which all the closed curves are homotopic to those in $\mathcal{P}$.
 This is not true for general hyperbolic metrics. One problem is that union of the geodesics obtained by replacing each curve $C_i$ of  $\mathcal{P}$ by its geodesic representative might not be closed, and there are others.

  In \cite{ALPSS} we discussed this question and we gave a necessary and sufficient condition on hyperbolic structures on surfaces of infinite type to have a pair of pants decomposition. One way of stating that result is to say that a hyperbolic metric satisfies this property if and only if it is the intrinsic metric of some conformal structure. This is the reason why in what follows we shall consider only hyperbolic metrics that are intrinsic in the sense we defined. A hyperbolic metric $S$ on $\Sigma$ has \emph{Fenchel-Nielsen coordinates} $\left((l_S(C_i),\theta_S(C_i))\right)_{i=1,2,\ldots}$ with respect to $\mathcal{P}$, using the notation of \cite{ALPSS}. For the convenience of the reader, the definition of these coordinates as well as the precise definitions of the three Teichm\"uller spaces that we mentioned above are recalled in \S \ref{sec:three} below.

We consider a conformal structure $S_0$ on $\Sigma$ which we call the basepoint of the Teichm\"uller space. We denote by $( \mathcal{T}_{qc}(S_0), d_{qc})$  the quasiconformal Teichm\"uller space equipped with the corresponding metric, and by $( \mathcal{T}_{ls}(S_0), d_{ls})$ the length-spectrum Teichm\"uller space equipped with its metric. We also let $\mathcal{P}=\{C_i\}_{i=1,2,\ldots}$ be a fixed pair of pants decomposition of $\Sigma$ and we denote by  $( \mathcal{T}_{FN}(S_0), d_{FN})$ the resulting Fenchel-Nielsen Teichm\"uller space equipped with its metric. The Fenchel-Nielsen Teichm\"uller space depends on the choice of $\mathcal{P}$, but we will not mark this dependence explicitly unless this is necessary. Hence the space $ \mathcal{T}_{FN}(S_0)$ and its metric are not intrinsic objects associated to $S_0$ but they constitute a useful tool to study the other spaces, because  $\mathcal{T}_{FN}(S_0)$  has explicit coordinates and it is isometric to the sequence space $\ell^{\infty}$. We shall recall the definitions in Section \ref{sec:three}. 

We note that in this paper we consider the \emph{reduced} Teichm\"uller space theory. This means that if the ideal boundary of $S_0$ is non-empty (see e.g. \cite{Nag} for the definition of the ideal boundary), a Teichm\"uller space of $\Sigma$ is a set of equivalence classes of marked Riemann surfaces up to homotopy, where the homotopy is free on the boundary components.

Given a hyperbolic structure $S$ and a simple closed curve $C$ on $\Sigma$, we denote by $l_S(C)$ the length of the unique $S$-geodesic in the homotopy class of $C$. In the case where $S$ is a conformal structure on $\Sigma$, then we denote by $l_S(C)$ the length of the unique geodesic in the homotopy class of $C$ with respect to the intrinsic metric associated to $S$.

  We say that a conformal structure $S$ is \emph{upper-bounded with respect to $\mathcal{P}$} if there exists a constant $M>0$ such that for any simple closed curve $C_i$ in $\mathcal{P}$, we have $l_S(C_i)\leq M$.

  We say that a conformal structure is \emph{upper-bounded} if it is upper-bounded with respect to some pair of pants decomposition, or if it is upper-bounded with respect to a pair of pants decomposition $\mathcal{P}$ which is understood.

  A  {\it marked conformal structure} (respectively a  {\it marked hyperbolic structure}) on $\Sigma$ is a pair $(f,S)$ where $S$ is a surface homeomorphic to $\Sigma$ equipped with a conformal (respectively a hyperbolic structure) and  $f:\Sigma \to S$ a homeomorphism. A marked conformal (respectively hyperbolic) structure on $S$ induces a conformal (respectively hyperbolic) structure on the surface $\Sigma$ itself by pull-back. Conversely, a conformal (respectively hyperbolic) structure on $S$ can be considered as a marked hyperbolic surface, by taking the marking to be the identity homeomorphism of $\Sigma$. using this formalism, an element of Teichm\"uller space is then an equivalence class of marked hyperbolic structures $(f,S)$ where the equivalence relation $\sim$ defined by $(f,S)\sim (f',S')$ if there exists an isometry $h:S\to S'$ such that $h\circ f$ is homotopic to $f'$.
  We shall use the notation $[f,S]$ to denote the equivalence class of the marked surface $(f,S)$.

In \cite[Theorem 8.10]{ALPSS}, we proved the following:

\begin{theorem}      \label{th:bilipschitz}
Let $S_0$ be a conformal structure on $\Sigma$, and suppose that $S_0$ is upper-bounded.  Then we have a set-theoretic equality $\mathcal{T}_{qc}(S_0)= \mathcal{T}_{FN}(S_0)$. Furthermore,
the identity map
\[j : (\mathcal{T}_{qc}(S_0), d_{qc}) \ni [f,S]\mapsto
\left((l_S(C_i),\theta_S(C_i))\right)_{i=1,2,\ldots}
 \in (\mathcal{T}_{FN}(S_0),d_{FN})\]
is a locally bi-Lipschitz homeomorphism.
\end{theorem}

Since the metric $d_{FN}$ on  the Fenchel-Nielsen Teichm\"uller space $\mathcal{T}_{FN}(S_0)$ makes this space isometric to the sequence space $\ell^\infty$, Theorem \ref{th:bilipschitz}
gives a locally bi-Lipschitz homeomorphism between the quasiconformal Teichm\"uller space $(\mathcal{T}_{qc}(S_0),d_{ls})$ and $\ell^{\infty}$. An analogous result was proved by Fletcher in \cite{Fletcher}, in the setting of non-reduced Teichm\"uller spaces.

Our main goal in this paper is to give a local comparison result between the Fenchel-Nielsen metric and  the length spectrum metric. The latter metric, in the setting of surfaces of infinite type, was first studied by Shiga in \cite{Shiga}.
A famous lemma due to Wolpert (see the exposition in \cite{Abikoff}) implies that for any hyperbolic surface $S_0$, we have a natural inclusion
\begin{equation}\label{eq:inclusion2}
\mathcal{T}_{qc}(S_0)\hookrightarrow \mathcal{T}_{ls}(S_0)
\end{equation}
given by the identity map, and that this map is $1$-Lipschitz, that is, for any two elements $S$ and $S'$ in $\mathcal{T}_{qc}(S_0)$, we have $d_{ls}(S,S')\leq d_{qc}(S,S')$.
We note by the way that in general, this inclusion map is not surjective (see \cite{LP} for an example).

  Theorem \ref{th:bilipschitz}, combined with Wolpert's result, gives the following:

\begin{theorem}   \label{th:biliandwolpert}
Let $S_0$ be a conformal structure on $\Sigma$ which is upper-bounded. Then, for any $S$ in $\mathcal{T}_{FN}(S_0)$,
there exists a neighborhood $N$ of $S$ in
$\mathcal{T}_{FN}(S_0)$ and a constant $C>0$ that depends only on $N$ such that for any $S'$ and $S''$ in $N$, we have
\[d_{ls}(S',S'')\leq C d_{FN}(S',S'').\]

\end{theorem}

Besides the upper-boundedness property for conformal structures, we shall use the following stronger property, which we call \emph{Shiga's property}, because it was used in a similar context in Shiga's paper \cite{Shiga}.

  We say that a conformal structure $S$ \emph{satisfies Shiga's property with respect to $\mathcal{P}$} if there exist two positive constants  $\delta$ and $M$ such that the following holds
  \begin{equation} \label{eqn:Shiga}
\forall C_i\in \mathcal{P},  \ \delta \leq l_S(C_i)\leq   M.
  \end{equation}

 Like for the upper-boundedness condition,  we shall say that a conformal structure  \emph{satisfies Shiga's property} if it  satisfies such a property for some pair of pants decomposition, or if it satisfies it  for a pair of pants decomposition which is understood.
 
The main result of this paper is the following, which proof appears at the end of Section \ref{s:FN} (Theorem \ref{Th:comparison}).

\begin{theorem}\label{Theorem:comparison}
Let $S_0$ be a conformal structure on $\Sigma$ satisfying Shiga's condition (\ref{eqn:Shiga}) and let $\mathcal{T}(S_0)_{qc}$ be the corresponding Teichm\"uller space. For any element $S$ of $\mathcal{T}_{qc}(S_0)$ and for any positive constant $D$, there exists a
positive real number $C$ that depends only on $\delta, M, D$ and $d_{ls}(S_0,S)$ such that if  and if two elements
 $S_1$ and $S_2$ of $\mathcal{T}_{qc}(S_0)$ are in the open ball of centre $S$ and radius $D$, then $d_{FN}(S_1,S_2)<Cd_{ls}(S_1,S_2)$.
\end{theorem}

From Theorems \ref{th:bilipschitz}, \ref{th:biliandwolpert} and \ref{Theorem:comparison}, we deduce the following.

\begin{theorem}\label{th:local-Shiga}
Let $S_0$ be a conformal structure on $\Sigma$ satisfying Shiga's condition.
Then we have a set-theoretic equality $\mathcal{T}_{qc}(S_0)= \mathcal{T}_{ls}(S_0)=\mathcal{T}_{FN}(S_0)$, and
the identity map between any two of the three spaces with their respective metrics $d_{qc},d_{ls}$ and $d_{FN}$ is locally bi-Lipschitz.
\end{theorem}

This implies in particular that under Shiga's condition, $\mathcal{T}_{ls}(S_0)$ is locally bi-Lipschitz equivalent to the sequence space $\ell^{\infty}$. It also implies that the Fenchel-Nielsen Teichm\"uller space $\mathcal{T}_{FN}(S_0)$, as a set, does not depend on the choice of the  pair of pants decomposition of $S_0$, and that the identity map between two Fenchel-Nielsen spaces with the same basepoint and corresponding to different pairs of pants decompositions is a bi-Lipschitz homeomorphism. (In particular, the topologies induced are the same.)

We also show below (Theorem \ref{th:not-bL} in Section \ref{s:ls}) that there exists a conformal structure $S_0$ on $\Sigma$ that does not satisfy Shiga's condition and such that the inclusion map between the quasiconformal Teichm\"uller space  $(\mathcal{T}_{qc}(S_0), d_{qc})$ and the length spectrum Teichm\"uller space $(\mathcal{T}_{ls}(S_0), d_{ls})$ is not locally bi-Lipschitz onto its image.
(Recall that by Wolpert's inequality there is always a set-theoretic inclusion  $\mathcal{T}_{qc}(S_0)\subset \mathcal{T}_{ls}(S_0)$.)

The above results and their proofs, although they are formulated for surfaces of infinite topological type, apply with little changes  to surfaces of finite topological type.
 In the latter case, all Teichm\"uller spaces coincide setwise. Some of the results we obtain here for surfaces of infinite type are known to be true for surfaces of finite type, but we also obtain some new results. We consider the case of surfaces of finite topological type in section \ref{s:finite}.

In the case of a surface of finite type with punctures and nonempty boundary, we studied in the papers \cite{LPST1} and \cite{LPST2} a metric on Teichm\"uller space which we called the \emph{arc metric} and which we denote by $\delta_L$, whose definition is analogous to the length spectrum metric, but which uses lengths of arcs instead of lengths of closed curves.
 For such a surface, we let $\mathcal{D}$ be the set of boundary components of $S$.

 We recall that for a surface of finite type and for $\epsilon >0$, the 
 $\epsilon$-thick part of Teichm\"{u}ller space, denoted by $\mathcal{T}_{\epsilon}$, is defined as the space
\[\mathcal{T}_{\epsilon}(S)= \{X\in \mathcal{T}(S) \ \vert \ \forall \gamma\in\mathcal{S}, l_X(\gamma) \geq \epsilon \}
\]
(Theorem \ref{th:arc2} below).

 For surfaces of finite type with nonempty boundary, for
 $\epsilon >0$ and $L\geq \epsilon$, we define the {\it $\epsilon_{0}$-relative} $\epsilon$-thick part of Teichm\"{u}ller space, denoted by $\mathcal{T}_{\epsilon, \epsilon_{0}}$, as the subset of the $\epsilon$-thick part of Teichm\"{u}ller space defined as
\[\mathcal{T}_{\epsilon,\epsilon_{0}}(S)= \{X\in \mathcal{T}(S) \ \vert \ \forall \gamma\in\mathcal{S}, l_X(\gamma) \geq \epsilon \hbox{ and }
\forall \delta\in\mathcal{D}, l_X(\gamma) \leq \epsilon_{0} \}
.\]
 We prove the following (Theorem \ref{th:arc} below):

\begin{theorem} Let $S$ be a topological finite type surface. For any $\epsilon>0$ and any $\epsilon_{0}>0$, the identity map between any two of the three metrics $d_{ls}$, $d_{qc}$ and $\delta_L$ on $\mathcal{T}_{\epsilon,\epsilon_{0}}(S)$ is globally bi-Lipschitz.
\end{theorem}

 In the case where the surface $S$ is of finite type and with empty boundary, then we have a similar statement for the
 $\epsilon$-thick part of Teichm\"{u}ller space (Theorem \ref{th:arc2} below).

\section{The three Teichm\"uller spaces}\label{sec:three}

In oder to make the paper self-contained and for the convenience of the reader, we recall the precise definitions of the three Teichm\"uller spaces that we associate to a surface of infinite topological type, namely, the  quasiconformal Teichm\"uller space, the Fenchel-Nielsen Teichm\"uller space and the length-spectrum Teichm\"uller space. These spaces were considered in the papers \cite{ALPSS}, \cite{completeness}, and \cite{LP}.

We start with the quasiconformal Teichm\"uller space $\mathcal{T}_{qc}(S_0)$. In this definition the hyperbolic metrics do not play a significant role, that is, we only use the
underlying Riemann surface structure of such a metric. More precisely,
the elements of  $\mathcal{T}_{qc}(S_0)$ are the homotopy classes of conformal structures $S$ on $\Sigma$ such that the identity map between $\Sigma$ equipped with $S_0$ and $S$ on the domain and on the target respectively  is homotopic to a quasiconformal map.
         The space $\mathcal{T}_{qc}(S_0)$ is equipped with the  \emph{quasiconformal metric}, also called the \emph{Teichm\"uller metric},
 in which for any two homotopy classes of conformal structures $(\Sigma,S)$ and $(\Sigma,S')$,  their  {\it quasiconformal distance} $d_{qc}(S,S')$ is defined as
        \begin{equation}\label{eq:qc}
        d_{qc}(R,R')=\frac{1}{2}\log \inf \{K(f)\}
        \end{equation}
        where  the infimum is taken over the set of quasiconformal dilatations  $K(f)$ of  quasiconformal  homeomorphisms $f:(\Sigma,S)\to (\Sigma,S')$ which are homotopic to the identity.
 Here, we are using the notation $(\Sigma,S)$ to say that $S$ is a structure (conformal or hyperbolic) on the surface $S$, with the marking being the identity map.

The conformal structure $S_0$ is the {\it basepoint} of $\mathcal{T}_{qc}(S_0)$.

We now recall the definition of the Fenchel-Nielsen Teichm\"uller spaces $\mathcal{T}_{FN}(S_0)$. In this definition we use the intrinsic hyperbolic metric associated to a conformal structure, see the discussion in the introduction regarding the pair of pants decomposition rendered geodesic with respect to the intrinsic hyperbolic metric. The definition of $\mathcal{T}_{FN}(S_0)$ is relative to the choice of a (topological) pair of pants decomposition
 $\mathcal{P}=\{C_i\}$ of $\Sigma$, and to the Fenchel-Nielsen coordinates associated to that decomposition.  The definition of the Fenchel-Nielsen parameters is similar to the one that is done in the case of surfaces of finite type, and we considered them in detail for surfaces of infinite type in \cite{ALPSS}.

Let $S$ be a (homotopy class of conformal) structure on $\Sigma$. To each homotopy class of closed geodesics $C_i\in \mathcal{P}$, we consider its {\it length parameter} $l_S(C_i)$ as defined in \S \ref{intro} above, and its {\it twist parameter} $\theta_S(C_i)$, which is  defined only if $C_i$ is not the homotopy class of a boundary component of $\Sigma$, as a measure of the relative twist amount along the geodesic in the class $C_i$ between the two generalized pairs of pants that have this geodesic in common. The twist amount per unit time along the (geodesic in the class) $C_i$ is chosen so that a complete positive Dehn twist along $C_i$ changes the twist parameter by addition of $2\pi$.

   Thus, for any conformal structure on $S$, its {\it Fenchel-Nielsen parameter} relative to $\mathcal{P}$ is the collection of pairs
\[\left((l_S(C_i),\theta_S(C_i))\right)_{i=1,2,\ldots}\]
where it
 is understood that if $C_i$ is homotopic to a boundary component,
 then there is no twist parameter associated to it, and instead of a pair
  $(l_S(C_i),\theta_S(C_i))$, we have a single parameter
$l_S(C_i)$.

Now given two conformal structures $S$ and $S'$ on $\Sigma$, their {\it Fenchel-Nielsen distance} (with respect to $\mathcal{P}$) is

\begin{equation}\label{def:FND}
{d_{FN}(S,S')=\sup_{i=1,2,\ldots} \max\left(\left\vert \log \frac{l_S(C_i)}{l_{S'}(C_i)}\right\vert, \vert l_S(C_i)\theta_S(C_i)-l_{S'}(C_i)\theta_{S'}(C_i)\vert \right) },
\end{equation}
again with the convention that if $C_i$ is the homotopy class of a boundary component of $\Sigma$, then there is no twist parameter to be considered.

Two conformal structures $S$ and $S'$ on $\Sigma$ are said to be {\it Fenchel-Nielsen bounded} (relatively to $\mathcal{P}$) if their Fenchel-Nielsen distance is finite. Fenchel-Nielsen boundedness is an equivalence relation.

   We say that two hyperbolic structures $S$ and $S'$ on $\Sigma$
  are equivalent if there exists an isometry $(\Sigma,S) \to (\Sigma,S')$ which is homotopic to the identity.
 Now given our basepoint $S_0$ of Teichm\"uller space, the {\it Fenchel-Nielsen Teich\-m\"uller space} with respect to $\mathcal{P}$ and with basepoint $S_0$, denoted by $\mathcal{T}_{FN}(S_0)$,  is the space of  equivalence classes of conformal structures that are Fenchel-Nielsen bounded from $S_0$ relative to $\mathcal{P}$.

 The function $d_{FN}$ defined above is a distance function on $\mathcal{T}_{FN}(S_0)$ and we call it the {\it Fenchel-Nielsen distance} relative to the pair of pants decomposition $\mathcal{P}$.
The map
\[
\mathcal{T}_{FN}(S_0) \ni H \mapsto {(\log(l_H(C_i))- \log(l_{S_{0}}(C_i))
, l_H(C_i)\theta_H(C_i))}_{i =1,2,\ldots}\in \ell^\infty \]
 is an isometric bijection between $\mathcal{T}_{FN}(S_0)$ and the sequence space $\ell^\infty$. It follows from known properties of $\ell^{\infty}$-norms that the Fenchel-Nielsen distance on $\mathcal{T}_{FN}(S_0)$ is complete.

  Finally, we recall the definition of the length-spectrum Teichm\"uller space  $\mathcal{T}_{ls}(S_0)$ with basepoint $S_0$. Again, in this definition we use the intrinsic hyperbolic metric associated to a conformal structure, see the discussion in the introduction.

We let $\mathcal{S}$ denote the set of homotopy classes of simple closed curves on $\Sigma$ that are not homotopic to a point or to a puncture (but they can be homotopic to a boundary component).
       We first define the {\it length-spectrum constant} $L(f)$ of a  homeomorphism  $f:(\Sigma,S_0)\to (\Sigma,S)$ where $S$ and $S'$ are two (homotopy classes of) conformal structures on $\Sigma$ as

\[
L(f)= \sup_{\alpha\in\mathcal{S}} \left\{ \frac{l_{S'}(f(\alpha))}{l_{S}(\alpha)},\frac{l_{S}(\alpha)}{l_{S'}(f(\alpha))}\right\}.
\]

 This quantity $L(f)$ depends only on the homotopy class of $f$, and we  say that $f$ is {\it length-spectrum bounded} if $L(f)<\infty$.

 We consider that two hyperbolic metrics $(\Sigma,S)$ and $(\Sigma,S')$ on $\Sigma$ are equivalent if there exists an isometry (or, equivalently, a length spectrum preserving homeomorphism) from $(\Sigma,S)$ to $(\Sigma,S')$ which is homotopic to the identity.
The {\it length spectrum Teichm\"uller space} $\mathcal{T}_{ls}(S_0)$ of $\Sigma$ with basepoint $S_0$ is the set of
homotopy classes of conformal structures $S$ on $\Sigma$ such that the identity map $\mathrm{Id}: (\Sigma,S_0)\to (\Sigma,S)$ is length-spectrum bounded.

 The {\it length-spectrum} metric $d_{ls}$ on $\mathcal{T}_{ls}(S_0)$ is defined by taking the distance $d_{ls}(S,S')$ between two points in that space to be
        \begin{equation}\label{eq:d-ls}
        d_{ls}(S,S')=\frac{1}{2}\log L(\mathrm{Id}).
        \end{equation}
        where $\mathrm{Id}$ is the identity map between $(\Sigma, S)$ and $(\Sigma, S')$ (we note that the length-spectrum constant of a length-spectrum bounded homeomorphism depends only on the homotopy class of such a homeomorphism.)

\section{On the Fenchel-Nielsen distance and the length spectrum distance}\label{s:FN}
Let $S$ be a hyperbolic structure on the surface of infinite topological type $S$ and let  $\mathcal{P}=\{C_i\}$ be a geodesic pair of pants decomposition of $S$.

\begin{lemma}\label{Lemma:angle}
Let $\delta <M$ be two positive constants such that each $C_i\in \mathcal{P}$ satisfies $\delta\leq l_S(C_i)\leq M$.
Then, for each $C_i\in \mathcal{P}$, we can find  a simple closed geodesic $\beta_i$ satisfying the following properties:

\begin{enumerate}
\item \label{t1} $\beta_i$ intersects $C_i$ in a minimal number of points (this number is one or two);
\item \label{t2} $\beta_i$ does not intersect $C_j$, for any $j\neq i$;
\item \label{t3}  there is a constant $L$ depending only on $\delta$ and $M$ such that
$l_S(\beta_i)<L$;
\item \label{t4} the sine of  the intersection angle (or of the two angles) of $\beta_i$ with $C_i$ is bounded from below by
a positive constant that depends only on $M$.
\end{enumerate}
\begin{proof}
Topologically, the curves $\beta_i$ are represented in Figure \ref{dual}.  \bigskip
  \begin{figure}[ht!]
\centering
\includegraphics[width=.40\linewidth]{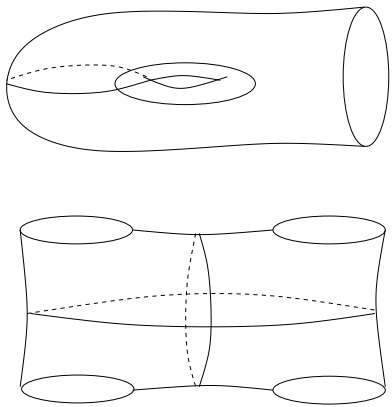}
\caption{\small {The curve $\beta_i$ used in the proof of Lemma \ref{Lemma:angle}.
In each case, we have represented the simple closed curves $C_i$ and $\beta_i$.}}
\label{dual}
\end{figure}
\bigskip
Using the inequalities $\delta\leq l_S(C_i)\leq M$, an
upper bound $L$ for $l_S(\beta_i)$ is obtained by estimates on hyperbolic right-angled hexagons and pentagons.
Using the upper bound  $L$ for $l_S(\beta_i)$ and the upper bound $M$ for $l_S(C_i)$, we can prove that $\sin\theta$
has a positive lower  bound depending on $L$ and $M$. We refer to Lemma 7.5 in \cite{ALPSS} for the details of the proof.

\end{proof}
\end{lemma}

Now we fix an element  $C_i\in \mathcal{P}$ and we let $\tau^t:S\rightarrow S_t$ be the time-$t$ Fenchel-Nielsen left twist deformation of $S$ along $C_i$. (At time $t$, we twist by an amount equal to $t$ measured on the curve $C_i$.)

Let $\beta$ be a simple closed geodesic on $S$. For all $t$ in $\mathbb{R}$, we denote by $\beta_t$ the simple closed geodesic  in $S_t$ homotopic to $\tau^t(\beta)$
and we let $l_t(\beta)= l_{S_{t}}(\beta_t)$ be its hyperbolic length. (Note that the class of $\beta_t$ is the same as the class of $\beta$ when we consider the hyperbolic structures as being on the same fixed base surface, or as marked surfaces with respect to a fixed base surface).
 The intersection
angle of $C_i$ and $\beta_t$ at a point $p\in C_i\cap \beta_t $
(measured from $C_i$ to $\beta_t$) is denoted by $\theta_t(p)$.

All angles used in this paper take their values in the  interval $[0,\pi]$.

We shall use the following formulae due to Wolpert \cite{Wolpert83}, concerning the first and second derivatives of the Fenchel-Nielsen flow. We use the formulation in  Weiss \cite{Weiss} p. 281).
\begin{lemma}\label{Lemma:Wolpert}
For any simple closed geodesic $\beta$, the function $t\mapsto l_t(\beta)$ is real-analytic, and we have
\[\frac{dl_t(\beta)}{dt}\big\vert_{t=0}=\sum\limits
_{p\in C_i \cap \beta}\cos \theta(p)\]
and
\[\frac{d^2l_t(\beta)}{dt^2}\big\vert_{t=0}=\sum_{p,q}
\frac{e^{l_1}+e^{l_2}}{2(e^{l(\beta)}-1)}\sin \theta(p) \sin \theta(q)
+
\sum_{p}
\frac{e^{l(\beta)}+1}{2(e^{l(\beta)}-1)}\sin^2 \theta(p).\]
 In the right hand side of the last inequality, the first sum is taken over the set of distinct ordered pairs of points $p,q$ in $C_i\cap \beta$, and  $l_1$ and $l_2$ are the lengths of the two subarcs that they subdivide on $\beta$,
and the second sum is taken over all points $p$ in $C_i\cap \beta$.
\end{lemma}

We shall use special cases of the above formulae, where $\beta$ intersects $C_i$ either in one or in two points.

In the case where $\beta$ and $C_i$ have only one intersection point $p$, with intersection angle $\theta (p)$,
the formulae become
$$\frac{dl_t(\beta)}{dt}\big\vert_{t=0}=\cos \theta(p)$$
and
$$\frac{d^2l_t(\beta)}{dt^2}\big\vert_{t=0}=
\frac{e^{l(\beta)}+1}{2(e^{l(\beta)}-1)}\sin^2 \theta(p).$$

In the case where $\beta$ and $C_i$ have two intersection points, denoted by $p_1$ and $p_2$, the formulae become
$$\frac{dl_t(\beta)}{dt}\big\vert_{t=0}=\cos \theta (p_1)+ \cos \theta (p_2)$$
and
\[\frac{d^2l_t(\beta)}{dt^2}\big\vert_{t=0}=
\frac{e^{l_1}+e^{l_2}}{(e^{l(\beta)}-1)}\sin \theta (p_1)\sin\theta(p_2)
+\frac{e^{l(\beta)}+1}{2(e^{l(\beta)}-1)}
(\sin^2 \theta(p_1)+ \sin^2 \theta(p_2)).\]

We now consider multi-twists, that is, composition of twist maps along a family of disjoint simple closed curves.

We let $t=(t_i)$, $i=1,2,\ldots$ be a sequence of real numbers, and for any real number $t$ we denote by $\tau^t:S\to S_t$ the multi-twist obtained by twisting an amount of $t_i$ along each curve $C_i$.

For $t=(t_i)$, $i=1,2,\ldots$, we set $\vert t\vert = \sup_{i=1,2,\ldots}\vert t_i\vert$.
\begin{proposition}\label{prop:twist}
Assume there exist two positive constants $\delta$ and $M$ such that each $C_i\in \mathcal{P}$ satisfies $\delta\leq l_S(C_i)\leq M$.
If $\displaystyle \sup_{\beta\in \mathcal{S}}|\log\frac{l_t(\beta)}{l(\beta)}|\leq D$,
then $$|t|<C\sup_{\beta\in \mathcal{S}}|\log\frac{l_t(\beta)}{l(\beta)}|,$$ where $C$ is a constant depending on $\delta, M$ and $D$.
\end{proposition}
\begin{proof}
It suffices to prove that for all $i=1,2,\ldots$, $\displaystyle\vert t_i\vert \leq C \sup_{\beta\in S}|\log\frac{l_t(\beta)}{l(\beta)}|$, where $C$ is a constant that depends on $\delta, M$ and $D$ and that does not depend on $i$.

For each $i$, we let $\beta_i$ be the simple closed geodesic given by Lemma $\ref{Lemma:angle}$. We shall apply the hypothesis $\displaystyle |\log\frac{l_t(\beta)}{l(\beta)}|\leq D$ to $\beta=\beta_i$ and show that
\begin{equation}\label{eq:for}\forall i=1,2,\ldots, \displaystyle\vert t_i\vert \leq C |\log\frac{l_t(\beta_i)}{l(\beta_i)}|.
\end{equation}

Note that the length $l_t(\beta_i)$ is affected by the twist along $C_i$, and not by any twist along $\beta_j$ for $j\not=i$..
 
 From \ref{eq:for}, we will then get 
 \[|t|=\sup\{|t_i|\} \leq C \sup_{\beta_i} |\log\frac{l_t(\beta_i)}{l(\beta_i)}|
 \leq \sup_{\beta\in \mathcal{S}}|\log\frac{l_t(\beta)}{l(\beta)}|,\]
 which is what we need to prove.
 
 Thus, we now prove \ref{eq:for}. We only need to assume that 
$$\sup_{\beta_i} |\log\frac{l_t(\beta_i)}{l(\beta_i)}|\leq D,$$
which is weaker than our assumption that 
 $$\sup_{\beta\in \mathcal{S}}|\log\frac{l_t(\beta)}{l(\beta)}|\leq D.$$

Without loss of generality, we can assume that $t_i>0$. In the following estimates, we can
restrict our attention to the pair of pants (or two pairs of pants) that contains $C_i$, 
since we only need to consider
the ratio $\frac{l_t(\beta)}{l(\beta)}$. We denote $t_i=t$ for simplicity.

There are two cases:

\noindent \textbf{Case I:}  $C_i$ intersects $\beta_i$ at a single point $p\in S$. Let $\theta$ be the angle at that intersection point.

By Lemma \ref{Lemma:angle}, there are positive constants $\rho_0=\rho_0(M)$  and $L=L(\delta,M)$ such that
$\sin\theta\geq \rho_0, l(\beta_i)<L$.
Since the function $t\mapsto l_t(\beta)$ is real-analytic, we can write
\[l_t(\beta)=l(\beta)+\frac{dl_t(\beta)}{dt}|_{t=0}t+\frac{d^2l_t(\beta)}{dt^2}|_{t=0}\frac{t^2}{2}+o(t^2).\]

From Lemma $\ref{Lemma:Wolpert}$, we obtain
\begin{equation}\label{Equation:W}
l_t(\beta_i)=l(\beta_i)+ \cos\theta \cdot t+ \frac{e^{l(\beta_i)}+1}{4(e^{l(\beta_i)}-1)}\sin^2\theta \cdot t^2+o(t^2).
\end{equation}

We now use the following result of Kerckhoff \cite{Kerckhoff}:

\begin{lemma} \label{lem:kerckhoff} The function $t\mapsto l_t(\beta_i)$ is strictly convex  and the function
$t\mapsto \cos\theta_t$ is strictly increasing.
\end{lemma}
In particular, if $l_t(\beta_i)$ attains its minimum at $t_0$, then
$\cos\theta_{t_0}=0$ (or, equivalently, $\theta_{t_0}=\frac{\pi}{2}$).
When $t<t_0$, $\cos\theta_t<0$ and when $t>t_0$, $\cos\theta_t>0$.

We set $\displaystyle |\log\frac{l_t(\beta_i)}{l(\beta_i)}|= \eta\leq D$. Then
$e^{-\eta}\leq\displaystyle  \frac{l_t(\beta_i)}{l(\beta_i)}\leq e^\eta$. Since
$1-e^\eta \leq e^{-\eta}-1\displaystyle \leq \frac{l_t(\beta_i)}{l(\beta_i)}-1\leq e^\eta -1$,
we have
\begin{eqnarray*}
|l_t(\beta_i)-l(\beta_i)|&=&l(\beta_i)|\frac{l_t(\beta_i)}{l(\beta_i)}-1| \\
&\leq& l(\beta_i)|e^{\eta}-1| \\
&=& l(\beta_i)(\eta+\sum_{n\geq 2} \frac{\eta^n}{n!}).
\end{eqnarray*}

By assumption, $l(\beta_i)<L$ and $\eta \leq D$. As a result,
\begin{equation}\label{equation:bound}
|l_t(\beta_i)-l(\beta_i)|\leq L(1+\sum_{n\geq 2} \frac{D^{n-1}}{n!})\eta= e(D)L \eta,
\end{equation}
where $e(D)$ is a constant that depends only on $D$.

Let $\lambda>0$ be a fixed sufficiently small positive constant, to be determined later.

First assume that  $\cos\theta\geq \lambda$. Applying the mean value theorem to the function $f(t)=l_t(\beta_i)$ on the interval $[0,t]$ and using the fact that  $f'(t)=\cos\theta_t$, we have
\begin{equation}\label{equation:mean}
l_t(\beta_i)-l(\beta_i)=\cos\theta_\xi t,
\end{equation}
for some $\xi\in [0,t]$.

Since $\cos\theta_\xi \geq  \cos\theta \geq \lambda$ (Lemma \ref{lem:kerckhoff}), combining $(\ref{equation:mean})$ with $(\ref{equation:bound})$, we have
$$t=\frac{l_t(\beta_i)-l(\beta_i)}{\cos\theta_\xi} \leq \frac{e(D)L\eta}{\lambda}.$$

If $|\cos\theta|< \lambda$, we let $\beta_i'$ be the unique geodesic on $S$ homotopic to the image of $\beta_i$ under the action of a positive Dehn twist along $C_i$.
Note that the hyperbolic
length of $\beta_i'$ is bounded by $L+M$ and, in fact, $\beta_i'= l_T(\beta_i)$, where $T=l_S(C_i)\geq \delta$. The value $T$ is the time needed for a full Dehn twist along $\beta_i$.

It is clear that $\beta_i'$ also satisfies the properties (1)-(3) in Lemma \ref{Lemma:angle}. Property (4) follows then from these three (see the proof of Lemma 7.5 in \cite{ALPSS}).
Let $p'$ be the intersection point of $C_i$ with $\beta_i'$ and $\theta'$ be the corresponding intersection angle.
Thus, there is a positive constant $\rho_1=\rho_1(L+M)$ such that
$\sin\theta'\geq \rho_1$. Let $\rho=\min \{\rho_0,\rho_1\}$.

We want to give a positive lower bound for $\cos\theta'$.

Since the hyperbolic length of $l_t(\beta_i), 0\leq t\leq T$ is bounded above by $L+M$ and since $\frac{e^x+1}{e^x-1}$ is a strictly decreasing function of $x$,
we have
\begin{equation}\label{equation:K}
\frac{e^{l_t(\beta_i)}+1}{e^{l_t(\beta_i)}-1}>\frac{e^{L+M}+1}{e^{L+M}-1},  \ \mathrm{for} \ 0\leq t\leq T.
\end{equation}

Let us set $\displaystyle K=\frac{e^{L+M}+1}{4(e^{L+M}-1)}$.

From Wolpert's formula (Lemma \ref{Lemma:Wolpert}), the second derivative with respect to $t$ of the length function  $l_t(\beta_i)$ is equal to
 \begin{equation}\label{equation:K1} \frac{e^{l_t(\beta_i)}+1}{2(e^{l_t(\beta_i)}-1)}\sin^2 \theta_t.
  \end{equation}

Inequality (\ref{equation:K}) shows that
   \begin{equation}\label{equation:K2}
  \forall t\in [0,T],  \ \frac{e^{l_t(\beta_i)}+1}{2(e^{l_t(\beta_i)}-1)}\sin^2\theta_t> K\sin^2\theta_t.
  \end{equation}

Thus, for $0\leq t\leq T$, the second derivative of $l_t(\beta_i)$ with respect to $t$ is bounded below by $ K\sin^2 \theta_t.$

For $0\leq t\leq T$,
we have $\sin\theta_t\geq \min \{\sin\theta, \sin\theta_T=\sin\theta' \}\geq \rho$, since $\sin \theta \geq \rho_0$ and $\sin \theta'\geq \rho_1$.

Thus, we have, using (\ref{equation:K2}),
$$\frac{d\cos\theta_t}{dt}=\frac{d^2l_t(\beta_i)}{dt^2}\geq K\rho^2.$$
As a result, and applying again the mean value theorem,
$$\cos\theta'-\cos\theta\geq K\rho^2T\geq K\rho^2\delta.$$
Now we set $\lambda=\frac{K\rho^2\delta}{2}$. Since
$|\cos\theta|<\lambda$ and $\cos\theta'-\cos\theta \geq 2\lambda$, we have $\cos\theta'>\lambda$.
The same arguments used in the first subcase show that

\[t\leq\frac{e(D)(L+M)\eta}{\lambda}.\]

The remaining subcase is when $\cos\theta\leq -\lambda$.

\begin{figure}[ht!]
\centering
\includegraphics[scale=0.65, bb=102 420 476 794]{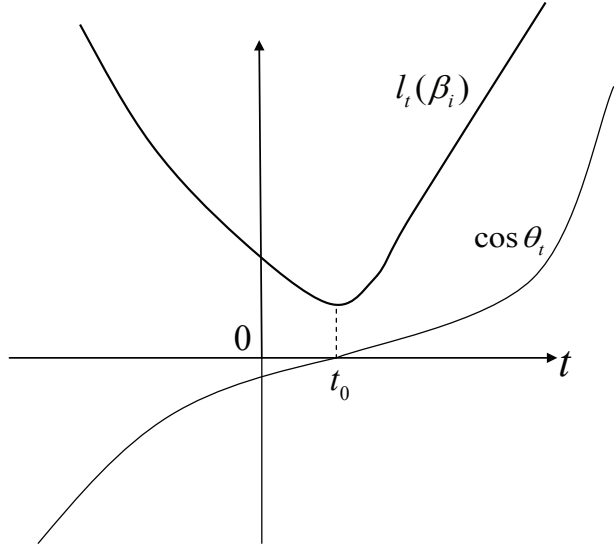}
\caption{\small{$l_t(\beta_i)$ is strictly convex and
$\cos\theta_t$ is strictly increasing.
If $l_t(\beta_i)$ attains its minimum at $t_0$, then
$\cos\theta_{t_0}=0$.}}
\label{fig:convex}
\end{figure}

Since $\sin\theta\geq \rho_0$, we  have $ -\sqrt{1-\rho_0^2} \leq \cos\theta <-\lambda$.
By Lemma \ref{lem:kerckhoff}, if $l_t(\beta_i)$ attains its minimum at $t_0$,
then $t_0>0$. This uses the fact that $\cos\theta <0$. See Figure \ref{fig:convex}. Set $N=[t_0]+1$. Since $l_t(\beta_i)$ decreases and
$\sin\theta_t$ increases when $t\in [0,t_0]$,
we have
$$\frac{e^{l_t(\beta_i)}+1}{e^{l_t(\beta_i)}-1}\sin\theta_t^2 \geq
\frac{e^{l(\beta_i)}+1}{e^{l(\beta_i)}-1}\sin\theta^2 \geq \frac{e^L+1}{e^L-1}\rho_0^2, \ \mathrm{for} \ t\in [0,t_0].$$
As a result, the first order derivative of $\cos\theta_t$ satisfies
$$\frac{d \cos\theta_t}{dt}\geq \frac{e^L+1}{4(e^L-1)}\rho_0^2, \ \mathrm{for} \ t\in [0,t_0].$$
Note that $t_0$ is exactly the value when $\cos\theta_t$ equals $0$, it follows that
$$\sqrt{1-\rho_0^2} \geq \cos\theta_{t_0}-\cos\theta \geq \frac{e^L+1}{4(e^L-1)}\rho_0^2 t_0.$$
This shows that  $N$ is bounded above by
$$\frac{4(e^L-1)}{(e^L+1)}\frac{\sqrt{1-\rho_0^2}}{\rho_0^2}+1.$$

Let $\beta_i^N$ be the geodesic on $S$ homotopic
to the image of of $\beta_i$ under $N$-order Dehn twist along $C_i$.
The intersection
angle $\theta_N$ of $C_i$ and $\beta_i^N$ satisfies $\cos\theta_N>0$
and the length of $\beta_i^N$ is bounded above by $L+NM$. By repeating the same argument as above, we complete the proof of
Case \textbf{(I)}.

\noindent \textbf{Case II}:   $\beta_i$ intersects $C_i$ at two different points $p_1, p_2$ (and we denote the intersecting angle by $\theta_1$ and $\theta_2$ respectively).

 In this case, we consider the formula
\begin{eqnarray*}\label{equation:main2}
l_t(\beta_i)&=&l(\beta_i)+ (\cos\theta_1+\cos\theta_2 )\cdot t \\
&+&\left( \frac{e^{l_1}+e^{l_2}}{(e^{l(\beta_i)}-1)}\sin \theta_1 \sin \theta_2 \cdot \frac{t^2}{2}
+ \frac{e^{l(\beta_i)}+1}{2(e^{l(\beta_i)}-1)}(\sin^2 \theta_1 + \sin^2 \theta_2)
\right)
+o(t^2).
\end{eqnarray*}
One checks that  $\displaystyle \frac{e^{l_1}+
e^{l_2}}{(e^{l(\beta_i)}-1)}\sin \theta_1 \sin \theta_2
\frac{e^{l(\beta_i)}+1}{2(e^{l(\beta_i)}-1)}(\sin^2 \theta_1 + \sin^2 \theta_2)
\geq A$, for some constant $A$ depending on $\rho_0$ and $L$.

Now fix a sufficiently small constant $0<\lambda_0< \frac{A}{2}\delta$. If $\cos\theta_1+\cos\theta_2 \geq \lambda_0$, using again the mean value theorem for the function $t\mapsto l_t(\beta_i)$,  it
 is easy to show that
 \[\vert t\vert < \frac{e(D)L\eta}{\lambda_0}.\]

If $|\cos\theta_1+\cos\theta_2 |< \lambda_0$, then we replace $\beta_i$ with its image $\beta'_i$ under the action of
positive Dehn twist along $C_i$. Let $\theta_1'$ and $\theta_2'$ be the intersection angles of $C_i$ and $\beta'_i$,
then the same proof as Case $\mathbf{I}$ shows that $|\cos\theta_1'+\cos\theta_2'|\geq \lambda_0$. As a result, we also have $|t|<\displaystyle  \frac{e(D)(L+2M)}{\lambda_0}\eta$. If $\cos\theta_1+\cos\theta_2 \leq -\lambda_0$, then we have to replace $\beta_i$
with the image of $\beta_i$ under the action of a $N-$order Dehn twist along $C_i$, denoted by $\beta_i^N$,
such that $l(\beta_i^N)<L+NM$ and the two intersection angles of $C_i$ and $\beta_i^N$ are non-negative. To give a upper bound
of $N$, we use the same argument as the last step of Case $(\mathbf{I})$, observing that
the two intersection angles have the same behavior under Fenchel-Nielsen twist deformation.
\end{proof}

Now we can prove the following

\begin{theorem}\label{Th:comparison}
Let $S_0$ be a conformal structure satisfying Shiga's condition
 (\ref{eqn:Shiga}), let $\mathcal{T}(S)$ be the corresponding Teichm\"uller space and let $S$ be a point in $\mathcal{T}_{ls}(S_0)$.
If $d_{ls}(S,S_1)\leq D$ and $d_{ls}(S,S_2) \leq D$ for some positive real number $D$, then $d_{FN}(S_1,S_2)<Cd_{ls}(S_1,S_2)$, where $C$ is a positive
constant that depends only on $\delta, M, D$.
\end{theorem}
\begin{proof}
There are positive constants $\delta_1$ and $M_1$, depending on $\delta, M$ and $D$, such that for each $C_i\in \mathcal{P}$, its hyperbolic length in $S_1$ satisfies
$$\delta_1\leq  l_{S_1}(C_i)\leq M_1.$$

By assumption, $d_{ls}(S_1,S_2)\leq 2D$. Then by Proposition \ref{prop:twist},
\begin{eqnarray*}
d_{FN}(S_1,S_2)&\leq & d_{FN}(S_1,S)+ d_{FN}(S,S_2) \\
&<& (C+1) d_{ls}(S_1,S_2)
\end{eqnarray*}
where $C$ is a positive constant depending on $\delta_1,M_1$ and $D$.
\end{proof}

 \section{The Teichm\"uller distance and the length spectrum distance}\label{s:ls}

   In this section, we show that the result in Theorem \ref{th:local-Shiga} is false if we remove Shiga's condition.

Let $S_0$ be a conformal structure on the surface (of infinite type) $\Sigma$ and $\mathcal{T}_{qc}(S_0)$
be its quasiconformal Teichm\"uller space, let $S$ be an element of $\mathcal{T}_{qc}(S_0)$ and
let $\alpha$ be a simple closed geodesic
on $S$. As before, we denote by $S_t$ be the hyperbolic surface obtained by by the time-$t$ Fenchel-Nielsen twist
deformation of $S$ along $\alpha$.
Recall the following proposition, which is a direct corollary of Lemma 7.4 in \cite{ALPSS}.

\begin{proposition}\label{prop:teich-lower} Let $T$ be a positive constant.
For $|t|<T$, we have
$$d_{qc}(S,S_t)\geq C |t|,$$
where $C$ is a positive constant depending only on $T$.
\end{proposition}

To compare the Teichm\"uller distance $d_{qc}(S,S_t)$ and the length spectrum distance
$d_{ls}(S,S_t)$, we show the following inequality.
\begin{lemma}\label{lem:ls-upper2}
$$d_{ls}(S,S_t) \le \frac{1}{2}\log
\sup_{\gamma \vert  i(\alpha, \gamma)\ne 0}\frac{i(\alpha,\gamma)|t|}{l_S(\gamma)}.$$
\end{lemma}
\begin{proof}
Without loss of generality, we can assume that $t>0$.
For any simple closed curve $\beta$ satisfying $i(\beta,\alpha)\not=0$, let $l_{S_t}(\beta)$
denote the hyperbolic length of $\beta$ in
$S_t$. From the definition of the Fenchel-Nielsen twist, we easily have
$$l_S(\beta)-i(\alpha,\beta)t \leq l_{S_t}(\beta)\leq l_S(\beta)+i(\alpha,\beta)t.$$

The length spectrum distance can be written as
$$d_{ls}(S,S_t) = \max \{\frac{1}{2}\log
\sup_\gamma\frac{l_{S_t}(\gamma)}{l_{S}(\gamma)},\frac{1}{2}\log
\sup_\gamma\frac{l_{S}(\gamma)}{l_{S_t}(\gamma)}\},$$
where the supremum is taken for all essential simple closed curves $\gamma$.

The hyperbolic length of a homotopy class of simple closed curves $\gamma$ satisfying $i(\alpha, \gamma)=0$ is invariant under the twist along $\alpha$.
As a result,  we have
$$d_{ls}(S,S_t)
= \max \{\frac{1}{2}\log
\sup_{\gamma \vert i(\alpha, \gamma)\ne 0}\frac{l_{S_t}(\gamma)}{l_{S}(\gamma)},\frac{1}{2}\log
\sup_{\gamma\vert i(\alpha, \gamma)\ne 0}\frac{l_{S}(\gamma)}{l_{S_t}(\gamma)}\}.
$$

For any simple closed curve $\gamma$ with $i(\alpha, \gamma)\ne 0$,
$$\log\frac{l_{S_t}(\gamma)}{l_{S}(\gamma)} \le|\log \frac{ l_S(\gamma)+i(\alpha,\gamma)t}{l_S(\gamma)}|\leq \frac{i(\alpha,\gamma)t}{l_S(\gamma)}$$
and
$$\log\frac{l_{S}(\gamma)}{l_{S_t}(\gamma)} \le|\log \frac{l_S(\gamma)} { l_S(\gamma)-i(\alpha,\gamma)t}|\leq \frac{i(\alpha,\gamma)t}{l_S(\gamma)}.$$
Then we have
$$d_{ls}(S,S_t) \le \frac{1}{2}
\sup_{\gamma\vert i(\alpha, \gamma)\ne 0}\frac{i(\alpha,\gamma)t}{l_S(\gamma)}.$$
\end{proof}
Note that if $l_S(\alpha)\leq L$, then it follows from the collar lemma that
there is a constant $C$ depending on $L$, such that for
any simple closed geodesic $\gamma$ with $i(\alpha, \gamma)\ne 0$,
$l_S(\gamma)\geq C i(\gamma,\alpha)|\log l_S(\alpha)|$.
Then Lemma \ref{lem:ls-upper2} gives:
\begin{lemma}\label{lem:ls-upper}
If $l_S(\alpha)\leq L$, then there is a constant $C$ depending on $L$ such that
$$d_{ls}(S,S_t) \le
\frac{t}{2C|\log l_S(\alpha)|}.$$
\end{lemma}
Combining Proposition \ref{prop:teich-lower} and Lemma \ref{lem:ls-upper}, we have
\begin{theorem}\label{thm:comparison}
If $l_S(\alpha)\leq L$ and $0<|t|<T$, then there exists a constant $C$ depending on
$L$ and $T$, such that
$$\frac{d_{qc}(S,S_t)}{d_{ls}(S,S_t)}\geq C |\log l_S(\alpha)|.$$
\end{theorem}

As an application, we show that
\begin{theorem}\label{th:not-bL}
If $S_0$ is a conformal surface of infinite topological type
with a pair of pants decomposition $\mathcal{P}=\{C_i\}$ such that there is a
subsequence of $\{C_{i_k}\}$ contained in the interior of $S_0$ whose
  hyperbolic lengths tend to zero, then the identity map between the
Teichm\"uller spaces $(\mathcal{T}_{qc}(S_0),d_{qc}) $ and its image in $(\mathcal{T}_{ls}(S_0), d_{qc})$ is not locally bi-Lipschitz.
\end{theorem}
\begin{proof}
By assumption, there is a subsequence $\{C_{i_k}\}$ of elements of $\mathcal{P}$ with hyperbolic length
$l_{S_0}(C_{i_k})=\epsilon_k\to 0$. For any fixed $t$, let $S_{k,t}$ be the hyperbolic
surface obtained by $t$ time Fenchel-Nielsen twist deformation of $S_0$ along $C_{i_k}$.
By Theorem \ref{thm:comparison}, there is a constant $C$, depending on the maximum of $\epsilon_k$ and $|t|$, such that
$$\frac{d_{qc}(S_0,S_{k,t})}{d_{ls}(S_0,S_{k,t})}\geq C |\log \epsilon_k|.$$
Since $\log \epsilon_k \to \infty$ as $\epsilon_{0}\to 0$, we have
$$\lim_{k\to\infty}\frac{d_{qc}(S_0,S_{k,t})}{d_{ls}(S_0,S_{k,t})}= \infty.$$

To see that the identity map between $d_{qc}$ and $d_{ls}$ is not locally bi-Lipschitz, we reason by
contradiction. Assume there are constants $C_1,C_2 $, such that for any
$S\in \mathcal{T}_{qc}(S_0)$, if $d_{ls}(S_0,S)\leq C_1$, then $d_{qc}(S_0,S)\leq C_2 d_{ls}(S_0,S)$.

Consider $S_{k,t}$ as above, and note that the Teichm\"uller distance is controlled by $t$. In fact,
if $|t|<T$ and $l_{S_0}(C_{i_k})\leq L$, we have
$$d_{qc}(S,S_{k,t})\leq C\vert t\vert,$$
where $C$ is a constant depending on $T$ and $L$. See \cite[Lemma 8.3]{ALPSS} for the proof.
As a result, for any $k$, we can choose $\vert t\vert $ sufficiently small such that $d_{ls}(S_0,S_{k,t})\leq C_1$.
However, we have shown that as $k\to \infty$,
$$\frac{d_{qc}(S_0,S_{k,t})}{d_{ls}(S_0,S_{k,t})}\to \infty,$$
which contradicts the assumption that $d_{qc}(S_0,S)\leq C_2 d_{ls}(S_0,S)$.

\end{proof}

The following is an analogous result, with a sequence $\{S_k\}$ in $\mathcal{T}_{qc}(S_0)$, such that
$d_{qc}(S_0,S_k)\to \infty$, while $d_{ls}(S_0,S_k)\to 0$.

\begin{example}
Let $S_0$ be a conformal structure of infinite type
with pants-decomposition $\mathcal{P}=\{C_i\}$, such that there is a
subsequence of $\{C_i\}$, contained in the interior of $S_0$,
with hyperbolic length $l_{S_0}(C_{i_k})=\epsilon_k$ tends to zero.
Let $S_k$ be the hyperbolic surface obtained by $t_k$ time Fenchel-Nielsen
twist deformation of $S_0$ along $C_{i_k}$. Here $\{t_k\}$ is sequence of  positive
constants tends to infinity and satisfying $\displaystyle \frac{t_k}{|\log \epsilon_k|} \to 0$.
Then it follows from the proof of Lemma \ref{lem:ls-upper2} that
$d_{ls}(S_0,S_k)\leq \displaystyle  \frac{t_k}{2C|\log \epsilon_k|}\to 0$.
On the other hand, the fact that $\displaystyle d_{qc}(S_0,S_k)\to \infty$ follows from
Lemma 7.2 in \cite{ALPSS}.

\end{example}

\section{The case of surfaces of finite type}\label{s:finite}

In this  section, we consider a hyperbolic surface $S=S_{g,m,n}$ of finite topological  type, of genus $g$ with $m$ punctures and $n$ boundary components. This means that when we equip such a surface with a conformal or a hyperbolic structure, then around each puncture, $S$ has a neighborhood that is conformally equivalent to a punctured disk, and around each boundary component, $S$ has a neighborhood which is conformally equivalent to an annulus. It is known that in this finite-type case we have the
 set-theoretic equalities $\mathcal{T}_{qc}(S)=\mathcal{T}_{ls}(S)=\mathcal{T}_{FN}(S)$, and we shall simply denote the Teichm\"uller space of $S$ by $\mathcal{T}(S)$ unless a  particular metric has to be specified.

 The reader will notice that Proposition \ref{prop:teich-lower}, Lemmas \ref{lem:ls-upper2} and \ref{lem:ls-upper} and Theorem \ref{thm:comparison} are valid for any Riemann surface, whether it has finite or infinite topological type.

 From Theorem \ref{thm:comparison}, we deduce the following
 \begin{corollary}
 For any non-elementary Riemann surface of finite type, the identity map between the  Teichm\"uller and the length spectrum metrics on $\mathcal{T}(S)$ is not a quasi-isometry.
 \end{corollary}

 This result, for surfaces of finite conformal type (that is, without boundary) was obtained independently and by other methods by Choi and Rafi in \cite{CR} and by Liu, Sun and Wei in \cite{LSW}. The result for surfaces of infinite topological type was obtained by Liu and Papadopoulos in \cite{LP}. The result for finite type surfaces with boundary is new.

 To state other results for surfaces of finite type, we recall the definition of a metric that we introduced in \cite{LPST1} on the Teichm\"uller space of such a surface. The definition of this metric uses the set of homotopic classes of arcs on $S$. Let us give the precise definition.

An {\it arc} in $S$ is the homeomorphic image of a closed interval whose interior is in the interior of $S$ and whose endpoints are on the boundary of $S$. All homotopies of arcs that we consider are relative to $\partial S$, that is, they leave the endpoints of arcs on the set $\partial S$.
An arc is said to be {\it essential} if it is not homotopic (relative to $\partial S$) to a map whose image is in $\partial S$.

We let $\mathcal{B}=\mathcal{B}(S)$ be the union of the set of homotopy classes
of essential arcs on $S$ with the set of homotopy classes of simple closed curves which are homotopic to boundary components.

Given an element $\gamma$ of $\mathcal{B}$ and an element $X$ of the Teichm\"uller space $\mathcal{T}(S)$, the \emph{length} of $\gamma$ with respect to $X$, denoted by $l_X(\gamma)$ is defined, in analogy with the length of an element of $\mathcal{S}$, as the length of the unique geodesic arc homotopic to $\gamma$ in a hyperbolic metric representing $X$.

 In \cite{LPST1} and \cite{LPST2} we studied the following metric on $\mathcal{T}(S_{g,m,n})$. For $X$ and $Y$ in this space, we set

\begin{equation}\label{eq:def5}
\delta_L(X,Y) = \log\max \left(\sup_{\gamma \in \mathcal{S}\cup\mathcal{B}}\frac{l_Y(\gamma)}{l_X(\gamma)},
\sup_{\gamma \in \mathcal {C}\cup \mathcal{B}}\frac{l_X(\gamma)}{l_Y(\gamma)}\right).
\end{equation}
 We showed that this function $\delta_L$ defines a metric, and that this metric is also given by
\begin{equation}\label{eq:def6}
 \delta_{L}(X,Y)=\log\max\left(
 \sup_{\gamma\in\mathcal{B}}\frac{l_{Y}(\gamma)}{l_{X}(\gamma)},
 \sup_{\gamma\in\mathcal{B}}\frac{l_{X}(\gamma)}{l_{Y}(\gamma)}
 \right).
\end{equation}

 We call $\delta_L$ the \emph{arc metric} on the Teichm\"uller space of the surface with boundary.

Any hyperbolic surface of finite type $S_{g,m,n}$ obviously satisfies Shiga's Condition (\ref{eqn:Shiga}), and Theorem \ref{Th:comparison} applies to such a surface.
Let $L$ be an upper bound for the hyperbolic length of the boundary geodesics of $S_{g,m,n}$. A result by Bers \cite{BB} shows that there exists a pants decomposition of $S$ with an upper bound $L_0$ for the lengths  of the decomposition curves, with $L_0$ depending only on $g,m,n$ and $L$.

 We shall use the following classical terminology.

Given a positive real number $\epsilon$, the $\epsilon$-thick part of the Teichm\"uller space of $S$, denoted by $\mathcal{T}_{\epsilon}(S)$, is defined as
\[\mathcal{T}_{\epsilon}(S)= \{X\in \mathcal{T}(S) \ \vert \ \forall \gamma\in\mathcal{S}, l_X(\gamma) \geq \epsilon  \}
.\]

We let $\mathcal{D}$ be the set of boundary components of $S$.
We shall use the following  terminology that we introduced in \cite{LPST2}.

For $\epsilon >0$ and $L\geq \epsilon$, the {\it $\epsilon_{0}$-relative} $\epsilon$-thick part of Teichm\"{u}ller space, denoted by $\mathcal{T}_{\epsilon, \epsilon_{0}}$, is the subset of the $\epsilon$-thick part of Teichm\"{u}ller space defined as
\[\mathcal{T}_{\epsilon,\epsilon_{0}}(S)= \{X\in \mathcal{T}(S) \ \vert \ \forall \gamma\in\mathcal{S}, l_X(\gamma) \geq \epsilon \hbox{ and }
\forall \delta\in\mathcal{D}, l_X(\delta) \leq \epsilon_{0} \}
.\]

We prove the following:

\begin{theorem}\label{th:arc} Let $S$ be a topological finite type surface. For any $\epsilon>0$ and any $\epsilon_{0}>0$, the identity map between any two of the three metrics $d_{ls}$, $d_{qc}$ and $\delta_L$ on $\mathcal{T}_{\epsilon,\epsilon_{0}}(S)$ is globally bi-Lipschitz.
\end{theorem}

\begin{proof} We first prove that the identity map \[\mathrm{Id}: (\mathcal{T}_{\epsilon,\epsilon_{0}},d_{ls})\to (\mathcal{T}_{\epsilon,\epsilon_{0}},d_{qc})\] is globally bi-Lipschitz.

It suffices to prove that for any $X,Y \in \mathcal{T}_{\epsilon, L}$, we have
\begin{equation}\label{pq}
d_{ls}(X,Y)\le d_{qc}(p,q)\le K d_{ls}(X,Y).
\end{equation}
where $K$ depending on the topological type of $S$, $\epsilon$ and $L$.

The left hand side inequality in (\ref{pq}) follows from Wolpert's lemma.

From Theorem \ref{th:local-Shiga} applied to surfaces of topological finite type, for any $D>0$, if $X,Y \in \mathcal{T}_{\epsilon, L}$ with $d_{ls}(x,y)\le D$, we have $d_{qc}(xX,Y)\le C d_{ls}(x,y)$, where $C$ depends on $D$, on the topological type of $S$ on $\epsilon$ and on $L$.
Therefore, if  $d_{ls}(X,Y)\le D$, the right hand side inequality of of (\ref{pq}) is satisfied. (We could take, for example, $D=1$.)

 Now assume that $d_{ls}(X,Y)\geq D$. From (12) of Theorem 6.3 in \cite{LPST1}, $d_{qc}(X,Y)\le \delta_{L}(X,Y)+D$.
From Theorem 3.6 in \cite{LPST2}, $\delta_L (X,Y) \le d_L(X,Y)+K$.
Thus, we have  $d_{qc}(X,Y)\le d_{ls}(X,Y)+D+K$.
This gives
$$
d_{qc}(X,Y) \le d_{ls}(X,Y)+K_1,
$$
where $K_1$ is a constant depending only on the topological type of $S$, on $\epsilon$ and on $L$.

As  $d_{ls}(X,Y)\geq D$, we have
$$
K_1= \frac{K_1} {D} D\le \frac{K_1} {D} d_{ls}(X,Y),
$$
and
$$
d_{qc}(X,Y) \le (1+\frac{K_1} {D})d_{ls}(X,Y).
$$
This proves the right hand side inequality of (\ref{pq}) in all cases.

It remains to show that the identity map \[\mathrm{Id}: (\mathcal{T}_{\epsilon,\epsilon_{0}},d_{ls})\to (\mathcal{T}_{\epsilon,\epsilon_{0}},\delta_L)\] is globally bi-Lipschitz. We use results proved in \cite{LPST1} and \cite{LPST2} on the natural embeddings between the Teichm\"uller space $\mathcal{T}(S)$ and the Teichm\"uller space $\mathcal{T}(S^d)$ of the double $S^d$ of $S$.
From the proof of Theorem 3.3 of \cite{LPST1}, this embedding is distance-preserving for the quasiconformal metrics on the two spaces. From Corollary 2.8 of \cite{LPST2}, this embedding is distance-preserving with respect to the metric $\delta_{L}$ on $\mathcal{T}(S)$ and  $d_{ls}$ on $\mathcal{T}(S^d)$. Furthermore, Proposition 4.2 of  \cite{LPST1} shows that the natural embedding
$\mathcal{T}(S)\to \mathcal{T}(S^d)$ sends an $\epsilon_{0}$-relative $\epsilon$-thick part of $\mathcal{T}(S)$ to an $\epsilon'$-thick part of $\mathcal{T}(S^d)$. We already showed that on such an $\epsilon'$-thick part of $\mathcal{T}(S^d)$, the identity map between the quasi-conformal and the length spectrum metrics is globally bi-Lipschitz. Therefore, the identity map between the quasi-conformal metric and the arc metric $\delta_L$ on the $\epsilon_{0}$-relative $\epsilon$-thick part of $\mathcal{T}(S)$ is globally bi-Lipschitz.

\end{proof}

In the case of a surface $S$ of finite type, we define more simply, for $\epsilon>0$,  the  $\epsilon$-thick part of $\mathcal{T}(S)$, denoted by $\mathcal{T}_{\epsilon}(S)$,  as
\[\mathcal{T}_{\epsilon}(S)= \{X\in \mathcal{T}(S) \ \vert \ \forall \gamma\in\mathcal{S}, l_X(\gamma) \geq \epsilon \}.
\]
We have the following theorem, analogous to Theorem \ref{th:arc}. The proof is similar to the proof of the first part of Theorem \ref{th:arc}.
\begin{theorem}\label{th:arc2} Let $S$ be a topological finite type surface without boundary. For any $\epsilon>0$, the identity map between the two metrics $d_{ls}$, $d_{qc}$ on $\mathcal{T}_{\epsilon}(S)$ is globally bi-Lipschitz.
\end{theorem}

 \end{document}